\magnification=1200
\hsize=16 true cm    \vsize=22.7 true cm
\hoffset=4 true mm   \voffset=2 true mm
\font\tengoth = eufm10
\font\sevengoth = eufm7
\font\fivegoth = eufm5
\newfam\gothfam
\textfont\gothfam=\tengoth
\scriptfont\gothfam=\sevengoth
\scriptscriptfont\gothfam=\fivegoth
\def\goth{\tengoth\fam\gothfam}
\font\tenmath =msbm10
\font\sevenmath=msbm7

\newfam\mathfam
\textfont\mathfam=\tenmath
\scriptfont\mathfam=\sevenmath
\scriptscriptfont\mathfam=\sevenmath
\def\math{\tenmath\fam\mathfam}
\def\build#1_#2{\mathrel{\mathop{\kern 0pt#1}\limits_{#2}}}
\def\K{{\math K}}
\def\wt{\hbox{wt}}
\def\dim{\hbox{dim}}
\def\wt{\hbox{wt}}
\def\m{{\bf m}}

\def\dim{\hbox{dim}}
\def\K{{\math K}}

\def\build#1_#2{\mathrel{\mathop{\kern 0pt#1}\limits_{#2}}}
\def\B{{\cal B}}\def\E{{\cal E}}

\def\qed{\hfill$\diamondsuit$}

\def\uqnm{U(\nm)}

\def\R{{\math R}}
\def\uq{U_q(\g)}\def\uqn{U_q(\n)}\def\uqnm{U_q(\nm)}\def\uo{U_q^0}\def\uqo{U
_q^0}\def\uqb{U_q(
\b)}\def\uqbm{U_q(\b^-)}
\def\lq{V_q}\def\lql{\lq(\lambda)}
\def\cqg{\C_q[G]}


\def\N{{\math Z}_{\geq 0}}\def\Z{{\math Z}}%
\def\n{{\goth n}}%
\def\C{{\math C}}%
\def\g{{\goth g}}%
\def\h{{\goth h}}%
\def\b{{\goth b}}\def\nm{{\goth n}^-}

\def\dim{\hbox{dim}}\def\tr{\hbox{Tr}}\def\wt{\hbox{wt}}

\def\P{{\cal P}}
\def\rp{\Delta^+}

\font\petittitre=cmbx8 scaled 1000 \font\grostitre=cmbx8 scaled 1440

\vskip 2cm\centerline{\grostitre{Adapted algebras for the
Berenstein-Zelevinsky conjecture}}
\centerline{\petittitre{ }}
\vskip 1cm\centerline{Philippe CALDERO}
\vskip 4cm\noindent
\par\noindent{\petittitre ABSTRACT.} {\it Let $G$ be a simply
connected semi-simple complex Lie group and fix
a maximal unipotent subgroup $U^-$ of $G$. Let $q$
be an indeterminate and let's denote by  $\B^*$ the
 dual canonical basis, [18], of the quantized algebra
 $\C_q[U^-]$ of regular functions on $U^-$. Following [19],
 fix a $\N^N$-parametrization of this basis, where $N=$\dim $U^-$.
 In [2], A. Berenstein
 and A. Zelevinsky conjecture that two elements of $\B^*$ $q$-commute
 if and only if they are multiplicative, i.e. their product
 is an element of $\B^*$ up to a power of $q$. For all reduced
 decomposition ${\tilde w_0}$ of the longest element of the Weyl group of
 $\g$, we associate a subalgebra $A_{\tilde w_0}$, called adapted
 algebra, of $\C_q[U^-]$ such
 that 1) $A_{\tilde w_0}$ is a $q$-polynomial algebra which
 equals $\C_q[U^-]$ up to localization, 2) $A_{\tilde w_0}$ is spanned by
 a subset of $\B^*$, 3) the Berenstein-Zelevinsky conjecture is
 true on $A_{\tilde w_0}$.
 Then, we test the conjecture when one element belongs to
 the $q$-center of $\C_q[U^-]$.}
  \vskip 2cm\noindent
This research has been partially supported by the EC TMR network
 "Algebraic Lie Representations", contract no. ERB FMTX-CT97-0100.
\vskip 1cm\noindent
\noindent
{\bf 0. Introduction.}\bigskip\noindent
{\bf 0.1.} Let $G$ be a simply connected semi-simple complex group, and let
$\g$ be its Lie algebra. Fix a Cartan subalgebra $\h$ of $\g$ and let
$\g=\n^-\oplus\h\oplus\n^+$ be the corresponding triangular decomposition.
For all $\lambda$ in the set of integral dominant weights $P^+$, let
$V(\lambda)$ be the simple $\g$-module with highest weight $\lambda$. \par
For all Lie algebra ${\goth a}$, the enveloping algebra of ${\goth a}$ will
be denoted by $U({\goth a})$. The module $V(\lambda)$ is naturally endowed
with a structure of $U(\g)$-module, and it is known that
$V(\lambda)=U(\n^-)v_\lambda$, where $v_\lambda$ is a highest weight vector
of $V(\lambda)$.
By the work of Lusztig, [19], there exists a base $\B$ of $U(\n^-)$ with the
following property : \smallskip\noindent
Set $\B(\lambda):=\{b\in \B,\,bv_\lambda\not=0\}$. Then, for all $\lambda$
in $P^+$, $\B(\lambda).v_\lambda$ is a base of $V(\lambda)$.
\smallskip\noindent
Let $U^-$ be the maximal unipotent subgroup of $G$ such that Lie($U^-)=\n^-$
and let $\C[U^-]$ be the algebra of regular functions on $U^-$. It is well
known
that $\C[U^-]$ is isomorphic to the symmetric algebra $S({\n^-}^*)$, and
thus to
$S(\n)$ via the Killing form. The left
multiplication of $U^-$ on itself provides a natural morphism $\phi$ from
$U(\n^-)$ on the algebra of differential operators on $\C[U^-]$. The pairing
on $U(\n^-)\times\C[U^-]$ such that $(u,f)$ is the value of $\phi(u)(f)$ on
the identity of $U^-$, is non degenerate.\par\noindent
Let $\B^*$ be the basis of $\C[U^-]$ dual to $\B$ for this pairing.
\smallskip\noindent
{\bf 0.2.} Berenstein and Zelevinsky have a conjecture in order to calculate
the basis $\B^*$. This conjecture involves the standard $q$-deformation of
the
basis $\B^*$. Let's present some notations.
The objects introduced in the previous section have a $q$-deformation. Let
$q$ be an indeterminate and let
$\uq$ be the quantized enveloping algebra on $\C(q)$ as defined by Drinfeld
and
Jimbo. Let $\uqn$ and $\uqnm$ be the quantized subalgebra corresponding
to $U(\n)$ and $U(\n^-)$, see 1.2, and let $\C_q[U^-]$ be the quantized
algebra of regular functions on $U^-$. For all $\lambda$ in $P^+$, let
$\lql$ be the simple $\uq$-module, 1.3, with highest weight vector
$v_\lambda$.
Then, there exists\smallskip\noindent
1) A canonical basis, still denoted by $\B$, of $\uqnm$, with the following
property\par\noindent
Set $\B(\lambda):=\{b\in \B,\,bv_\lambda\not=0\}$. Then, for all $\lambda$
in $P^+$, $\B(\lambda).v_\lambda$ is a base of $\lql$.\smallskip\noindent
2) A pairing on $\uqnm\times\C_q[U^-]$ which quantifies the pairing in 0.1
and which
permits to define the dual canonical basis $\B^*$ in
$\C_q[U^-]$.\smallskip\noindent
Remark that we can recover both basis  of 0.1 by specializing $q$ on 1.
Remark also that $\C_q[U^-]$ is isomorphic to $\uqn$.\par
We say that two elements $u$ and $u'$ in $\C_q[U^-]$ {\it $q$-commute} if
there exists
an integer $m$ such that $uu'=q^mu'u$. Two elements $b^*$ and $b'^*$ in
$\B^*$ are said to be multiplicative if there exists an $m$ such that
$q^m b^*b'^*\in\B^*$. We can formulate the Berenstein-Zelevinsky conjecture.
\medskip\noindent
{\bf Conjecture.} Two elements of the dual canonical basis $\B^*$
$q$-commute
if and only if they are multiplicative.\medskip\noindent
Let $A$ be a subalgebra of $\cqg$ such that $A$ is generated as a space by a
subset of
$\B^*$. In the sequel, we say that the Berenstein-Zelevinsky conjecture is
true
on $A$ if the assertion is true when both elements are in
$A$.
    \smallskip\noindent
{\bf 0.3.} We present results about the conjecture. Recall that the elements
of the
dual canonical basis are labelled by elements in $\N^N$, where $N$ is the
dimension
of $\n$. This is the so-called Lusztig parametrization of $\B^*$
which  depends on the choice of a
reduced decomposition ${\tilde w_0}$
of the longest element $w_0$ of the Weyl group.
 \par
In [21], the author uses the quiver approach of quantum groups.
First of all, he obtains the only if part of the Berenstein-Zelevinsky
conjecture,
[{\it loc. cit.}, Corollary 4.5], in the simply laced case.
Maybe the most remarkable result of [21] is obtained in the case A$_n$.
It can be interpreted as follows :
for all quiver orientation
$\Omega$ of the diagram of A$_n$, the author proves that there exists
 a subalgebra $A_\Omega$ of maximal GK-dimension
 of $\C_q[U^-]$ which is compatible with
$\B^*$ and on which the Berenstein-Zelevinsky conjecture is true. Moreover,
each pair of elements in $\B^*\cap A_\Omega$ $q$-commute (and so are
multiplicative)
and the previous parametrization yields a correspondance between
$\B^*\cap A_\Omega$ and the set ${\cal C}_\Omega$
of integral points of a simplicial cone in $\R^N$.\smallskip\noindent
Let $\g$ of type A$_n$. The non zero classical minors in $\C[U^-]$ belong
 to the (classical) dual canonical basis. They have a quantum analogue in
 $\C_q[U^-]$ which belong to the (quantum) dual canonical basis. A natural
question
 comes up : is the assertion of the Conjecture true when both elements are
 quantum minors. In [13], Leclerc, Nazarov and Thibon give a positive answer
 to this question. \smallskip\noindent
{\bf 0.5.} The main results of this article are\smallskip\noindent
1) a generalization of the result
of [21] for all semi-simple Lie algebra $\g$, see Theorem 2.2.
An adapted subalgebra of $\C_q[U^-]$, see Definition 2.2, is an algebra
which is spanned by a maximal subset $P^*$ of
$\B^*$ such that two elements in $P^*$ are multiplicative.
For all reduced decomposition $\tilde w_0$ of $w_0$, we provide an adapted
algebra $A_{\tilde w_0}$ which generalizes the algebra $A_\Omega$
discussed in 0.3. Recall that
in [5] and [6], we provide, for all reduced decomposition
of $w_0$, an subalgebra $A_{\tilde w_0}$ of  $\C_q[U^-]$ such that
a) the algebra $A_{\tilde w_0}$ is generated by a family $S_{\tilde w_0}$ of
$N$ algebraically independant elements which pairwise $q$-commute  b)
$A_{\tilde w_0}$ and  $\C_q[U^-]$ are equal up to localization by the
multiplicative part generated by $S_{\tilde w_0}$. Using Kashiwara's crystal
basis and some standard monomial theory, we prove that $A_{\tilde w_0}$ is
spanned by a subset of $\B^*$ and that the Conjecture is true on $A_{\tilde
w_0}$.
Then, we prove that $A_{\tilde w_0}$ is an adapted algebra.
\par\noindent
2) we test the Berenstein-Zelevinsky conjecture on the $q$-center of $\C_q[U^-]$, i.e. the space
generated
by elements which $q$-commute with all homogeneous elements of
$\C_q[U^-]$. Indeed, as the $q$-center of $\C_q[U^-]$ is spanned by a part
of
the dual canonical basis, a natural test for the conjecture can be stated as
follows :
is the conjecture true if one element belongs to the $q$-center. We give a
positive
answer to this question.
Moreover, it reduces the conjecture to a subset of $\B^*$, see
Corollary 3.2.
 \par\noindent
3) we study the case when $\g$ is of type B$_2$, see 3.4, and we remark that
in this case,
$\C_q[U^-]$ is a direct sum of adapted algebras.\medskip\noindent
To conclude, we give a presentation of these results in terms of cones in
$\N^N$, which are sets of parametrizations for the part $\P_{\tilde w_0}^*$
associated to $A_{\tilde w_0}$. These cones generalize the cones $\C_\Omega$
discussed in 0.4.     \bigskip\noindent
{\bf 1. Notations and recollection on global base.}\bigskip\noindent
{\bf 1.1.}  Let $\g$ be a semi-simple Lie $\C$-algebra  of rank $n$.
We fix a Cartan subalgebra   $\h$ of $\g$. Let $\g=\nm \oplus\h \oplus\n$
be the triangular decomposition  and let $\{\alpha_i\}_i$ be a base of the
root system
$\Delta$ resulting from this decomposition. Let $\rp$ be the set of positive
roots.
Let  $P$ be the weight lattice generated by the fundamental weights
$\varpi_i$, $1\leq i\leq n$.
Let $P^+:=\sum_i\,\N\varpi_i$ and $P:=\sum_i\,\Z\varpi_i$ be respectively
the semigroup of integral dominant weights and the group of integral
weights.
We endow $P^+$ with the ordering $\leq$ defined by
$\lambda\leq\mu\Leftrightarrow
\mu-\lambda\in P^+$.
 Let $W$ be the  Weyl group, generated by the
 reflections corresponding to the simple roots $s_i:=s_{\alpha_i}$, with
longest element
 $w_0$.
We note  $<\,,\,>$ the $W$-invariant form on $P$.
\medskip\noindent
{\bf 1.2.} Let $d$ be an integer such that $<P,P>\subset (2/d)\Z$.
Let $q$ be a indeterminate and set $\K=\C(q^{1/d})$.
Let $\uq$ be the simply connected quantized
enveloping  algebra on $\K$, as defined in [10]. Set
$d_i=<\alpha_i,\alpha_i>/2$
and $q_i=q^{d_i}$ for all $i$.
Let $\uqn$, resp. $\uqnm$,
be the subalgebra generated by the canonical generators
$E_{\alpha_i}$, resp. $F_{\alpha_i}$,
 of positive,  resp. negative, weights and the quantum Serre relations.
For all $\lambda$
in $P$, let $K_\lambda$ the corresponding  element in the
algebra $\uo=\K[P]$ of the torus of $\uq$.
 Recall the triangular decomposition
$\uq=\uqnm\otimes\uo\otimes\uqn$. We define the following subalgebras of
$\uq$ :
$$\uqb=\uqn\otimes\uo, \hskip 15mm \uqbm=\uqnm\otimes\uo.$$\par
$\uq$ is endowed with a  structure of Hopf algebra and the
comultiplication $\Delta$,
the antipode $S$ and the augmentation $\varepsilon$ are given by
$$\Delta E_i=E_i\otimes 1+K_{\alpha_i}\otimes E_i,\,\Delta F_i
=F_i\otimes K_{\alpha_i}^{-1}+ 1\otimes F_i, \Delta K_\lambda
=K_\lambda\otimes K_\lambda$$
$$S(E_i)=-K_{\alpha_i}^{-2}E_i,\,S(F_i)=-F_iK_{\alpha_i}^{2}, S(K_\lambda)
=K_{-\lambda}$$
$$\varepsilon(E_i)=\varepsilon(F_i)=0,\,\varepsilon(K_\lambda)
=1.$$\medskip\noindent Let's
denote by ad the adjoint action.
Remark that there exists an algebra isomorphism
$\uqn\simeq\uqnm$ which sends $E_{\alpha_i}$ on $F_{\alpha_i}$. \par
If $\alpha=\sum_im_i\alpha_i$ is in $Q^+:=\sum\N\alpha_i$, resp.
$Q^-:=\sum\Z_{\leq 0}\alpha_i$,
 then an element $X$ of the subspace of $\uqn$, resp. $\uqnm$, generated by
 $\{E_{i_1}^{n_1}\ldots E_{i_k}^{n_k},
 n_{i_1}\alpha_{i_1}+\ldots + n_{i_k}\alpha_{i_k}=\alpha\}$,
 resp. $\{F_{i_1}^{n_1}\ldots F_{i_k}^{n_k},
 n_{i_1}\alpha_{i_1}+\ldots + n_{i_k}\alpha_{i_k}=-\alpha\}$ will be
 called element of weight $\alpha$.
In this case, we set $\wt(X)=\alpha$ and $\tr(X)=\sum m_i$. \par
For each $u$ in $\uq$ we set $u_{(1)}\otimes
u_{(2)}=\Delta(u)
\in\uq\otimes\uq$. There exists a unique form, [22], [24],
$(\,,\,)$ on $\uqb\times\uqbm$
such that :
$$(E_i,\,F_i)=\delta_{ij}(1-q_i^2)^{-1},$$
$$(u^+,\,u_1^-u_2^-)=(\Delta(u^+),\,u_1^-\otimes u_2^-)\,,\hskip 1cm
u^+\in\uqb\,;\,u_1^-,u_2^-\in \uqbm$$
$$(u_1^+u_2^+,\,u^-)=(u_2^+\otimes u_1^+,\,\Delta(u^-))\,,\hskip 1cm
u^-\in\uqbm\,;\,u_1^+,u_2^+\in\uqb$$
$$(K_\lambda,\,K_\mu)=q^{-(\lambda,\mu)}\,,(K_\lambda,\,F_i)=0\,,
(E_i,\,K_\lambda)=0\,,\hskip 1cm
\lambda,\mu\in P$$
For all $\beta$ in $Q^+$, let $\uqn_\beta$, resp. $\uqnm_{-\beta}$, be the
subspace of $\uqn$, resp. $\uqnm$, with weight $\beta$, resp. $-\beta$.
The form $(\,,\, )$ is non degenerate on $\uqn_\beta
\times \uqnm_{-\beta}$, $\beta\in Q^+$. We have,  :
$$(XK_\lambda,\,YK_\mu)=q^{-(\lambda,\mu)}(X,\,Y)\,,
\hskip 1cm X\in \uqn,\,Y\in \uqnm$$
We can define a bilinear form $<\, ,\, >$ on
$\uq\times\uq$ by :
$$<X_1K_\lambda S(Y_1),\, Y_2K_\mu S(X_2)>=(X_1,\,Y_2)
(X_2,\,Y_1)q^{-(\lambda,\mu)/2}$$\noindent
where $X_1,X_2\in \uqn,\,Y_1,Y_2\in \uqnm
,\,\lambda,\mu\in P$.
Moreover, this form is non dege\-ne\-ra\-te.\par
We define the ring automorphism $x\mapsto \eta (x)$ of $\uqn$, resp.
$\uqnm$, such that
$\eta(q)=q^{-1}$ and  $\eta(E_i)=E_i$, resp. $\eta(F_i)=F_i$
for all $i$,
$1\leq i\leq n$. We also define the antihomomorphism
$\sigma$ of the $\C(q)$-algebra $\uqn$, resp. $\uqnm$,
 such that $\sigma(E_i)=E_i$, resp. $\sigma(F_i)=F_i$, for all
$i$,
$1\leq i\leq n$.
\smallskip\noindent
{\bf 1.3.} For all $\uqo$-module $M$, set $M_\lambda
:=\{m\in M,\, K_\mu.m=q^{(\mu,\lambda)}m\}$. Elements of $M_\mu$ are called
elements of weight
$\mu$. For all $\lambda$ in $P^+$, $\lql$ denotes the simple $\uq$-module
with
highest weight $\lambda$, with highest weight vector $v_\lambda$. It is
known, see [10, 4.3.6], that this module verifies the Weyl character formula.
In particular, for all $w$ in $W$,
dim$\lql_{w\lambda}=1$, and we can define extremal vectors $v_{w\lambda}$
in $\lql$ as in the classical case.\par
The dual space $V(\lambda)^*$ is endowed with a structure of left
$\uq$-module
by twisting with the antipode.
\medskip\noindent
{\bf 1.4.}  Fix a reduced decomposition $\tilde w_0=s_{i_1}\ldots s_{i_N}$
of $w_0$. Set $\beta_s:=s_{i_1}
\ldots s_{i_{s-1}}(\alpha_{i_s})$, $1\leq s\leq N$. Recall that $\{
\beta_s,\, 1\leq s\leq N\}$ is exactly the set of positive roots.
 We define
a total ordering of the positive roots as follows :
$$\beta_{1}>\beta_{2}>\ldots>\beta_{N}.$$
For $i$, $1\leq i\leq n$, let $T_i$ be the Lusztig automorphism, [19,
37.1.3], [23],
associated to $i$.
 We introduce the following elements in $\uqn$ :
 $E_{\beta_s}=T_{i_1}\ldots T_{i_{s-1}}(E_{i_s})$, $1\leq s\leq N$. Remark
that these elements
 depend on $\tilde w_0$. We index the so-called Poincar\'e-Birkoff-Witt
basis
 of $\uqn$ by $\m$ in $\N^N$
 as follows
 $$E(\m)=E_{\tilde w_0}(\m)=\prod_{i=N}^1{1\over [m_i]_{q_i}!}
 E_{\beta_i}^{m_i},$$
 where $[n]_{q_i}!=[n]_{q_i} [n-1]_{q_i}\ldots[1]_{q_i}$,
 $[n]_{q_i}={q_i^n-q_i^{-n}\over q_i-q_i^{-1}}$
 $q_i=q^{{(\beta_i,\beta_i)\over 2}}$.
We define in the same way a PBW-basis $(F(\m))_{\m\in \N^N}$
of $\uqnm$ via the isomorphism in
1.2.
 \smallskip\noindent
 Let $(E(\m)^*)_{\m\in \N^N}$ be the dual basis of the PBW-basis
$(F(\m))_{\m\in \N^N}$ for
$(\,,\,)$.
 By [14], we have
 $$E(\m)^*=\prod_{i=1}^N {\phi_{m_i}(q_i^2)\over (1-q_i^2)^{m_i}}E(\m),$$
 where  $\phi_m(z)=\prod_{k=1}^m (1-z^k).$
 \par\noindent
 We know by [9, Lemma 1.7] that the lexicographical
 ordering on the PBW-basis provides a filtration on $\uqn$ whose
 associated algebra is a $q$-polynomial algebra :\medskip\noindent
 {\bf Theorem.} {\it Fix a reduced decomposition $\tilde w_0$ of $w_0$ and
 set
 $${\cal F}_{\m}(\uqn)=\oplus_{{\bf n}\prec \m}\K E_{\tilde w_0}({\bf n}),
 \hskip 5mm
 \m\in\N^N,$$
 where $\prec$ is the lexicographical ordering of $\N^N$.
 Then, the associated graded algebra Gr$_{\tilde w_0}(\uqn)$ is generated
 by Gr$_{\tilde w_0}(E_\alpha)$, $\alpha\in\rp$ and with relations :
  $$\hbox{Gr}_{\tilde w_0}(E_\alpha)\hbox{Gr}_{\tilde w_0}(E_\beta)
  =q^{<\alpha,\beta>}\hbox{Gr}_{\tilde w_0}(E_\beta)\hbox{Gr}_{\tilde
w_0}(E_\alpha),
  \hskip 5 mm \alpha<\beta.$$
  In particular $\hbox{Gr}_{\tilde w_0}(E^*(\m))
  \hbox{Gr}_{\tilde w_0}(E^*({\bf n}))=\hbox{Gr}_{\tilde w_0}(E^*(\m+{\bf
n}))$
  up to a power of $q$. \qed}
 \medskip\noindent
 {\bf 1.5.} \def\L{{\cal L}}\def\B{{\cal B}}
 Let's define  now the following sub-$\Z[q]$-lattices of $\uqn$ :
 $$\L=\bigoplus_{\m\in\N^N}\Z[q]E(\m),\hskip 5 mm
\L^*=\bigoplus_{\m\in\N^N}\Z[q]E(\m)^*$$
 A remarkable result of Lusztig states that these lattices do not depend on
a
 reduced decomposition of $w_0$.
 The following theorem is due to Lusztig [17]. It introduces the
so-called
 canonical basis $\B$ of $\uqn$ which will be identified with the canonical
basis of
of $\uqnm$.\smallskip\noindent
 {\bf Theorem.} {\it Fix a reduced decomposition $\tilde w_0$ of $w_0$.
Then,
 for each $\m$ in $\N^N$,
 there exists a unique element $B(\m)=B_{\tilde w_0}(\m)$ in $\uqn$ such
that
 $\eta({B(\m)})=B(\m)$ and $B(\m)\in E(\m)+q\L$. The set $\B:=\{B(\m),\;
 \m\in \N^N\}$ is a basis of $\uqn$ which does not depend on $\tilde
w_0$.}\qed
 \smallskip\noindent
 Hence, for a fixed $\tilde w_0$, the bijection $\N^N\rightarrow \B$
 gives rise to a parametrization of $\B$. This will be called the Lusztig's
parametrization
 of the canonical basis $\B$.\par
 Let $\B^*=(B(\m)^*)_{\m\in\N^N}$ be the dual canonical basis of $\uqn$.
 As in [13, Proposition 16], we have\smallskip\noindent
 {\bf Proposition.} {\it Fix a reduced decomposition $\tilde w_0$ of $w_0$.
 Then, for each $\m$ in $\N^N$, the element $B(\m)^*$ is the unique element
$X$
 of $\uqn$ with weight $\wt(X)=\sum m_i\beta_i$ such that
 $$\eta(X)=(-1)^{tr(X)}q^{<\wt(X),\wt(X)>/2}q_{X}\sigma(X),\;\;
 X=E(\m)^*+q\L^*,$$
 where $q_X=\prod_iq_i^{m_i}$,
$\wt(X)=\sum_im_i\alpha_i$. }\medskip\noindent
{\bf Remark.} Suppose that $\tilde w_0$ corresponds to a quiver orientation.
Then, by [21, Lemma 4.1], see also [20, par. 7] for the non simply
laced case,
$$\hbox{Gr}_{\tilde w_0}(B^*(\m))=\hbox{Gr}_{\tilde w_0}(E^*(\m)).\leqno
(1.5.1)$$
\smallskip\noindent
Recall that two elements of the dual canonical basis are called
multiplicative if their product is an element of the dual canonical basis,
up
to
a power of $q$.
The following corollary of the proposition is due to Reineke, [21,
Corollary 4.5]. It
is an implication in the Berenstein-Zelevinsky conjecture. His proof is also available for the
non simply laced case.\medskip\noindent
{\bf Corollary.} {\it If two elements of the dual canonical basis are
multiplicative, then they $q$-commute.}\qed
 \medskip\noindent
 {\bf 1.6.} In this section we consider the canonical basis $\B$ in $\uqnm$.

Let $\tilde E_i$, $\tilde F_i$ : $\uqnm\rightarrow\uqnm$ be the
Kashiwara
operators, [19]. For $b\in\B$, $\tilde E_i(b)$, resp.
$\tilde F_i(b)$, equals some $b'$ in $\B\cup\{0\}$ modulo
$q\Z[q]\B$.
The rule $b\mapsto b'$ defines maps  $\tilde e_i$ and $\tilde f_i$ from
$\B$ to
$\B\cup\{0\}$.
For $b\in\B$, $1\leq i\leq n$, set $\varepsilon_i(b)=\hbox{Max}\{r,\,
\tilde e_i^r(b)\not=0\}$, and
$\E(b)=\sum_{i=1}^n\varepsilon_i(b)\varpi_i$.
\par Fix $\lambda$ in $P^+$ and define the following subset of $\B$
$$\B(\lambda)=\{b\in \B,\,\E(\sigma(b))\leq\lambda\}.$$
Let $L(\lambda)$ be the $\C[q]$-sublattice of $\lql$ generated by
$\{bv_\lambda,\,b\in \B(\lambda)\}$ and define the following subset
of $L(\lambda)/qL(\lambda)$
$$B(\lambda)=\{b \hbox{ mod } qL(\lambda),\, b\in \B(\lambda)\}.$$
In the sequel, the element $b$ mod $qL(\lambda)$ will be denoted
$[b]_\lambda$ or $[b]$ if no confusion occurs.
The following is well
known, [11, 12.3] :\medskip\noindent
{\bf Theorem.} {\it For all $\lambda\in P^+$, we have :\par\noindent
\item {(i)} the set $\B(\lambda)$ is the set of
elements $b$ of $\B$ such that $bv_\lambda$ is non zero
in $\lql$,
moreover $\B(\lambda)v_\lambda$ is a basis of $\lql$,\par\noindent
\item {(ii)} the pair $(L(\lambda), B(\lambda))$ is a crystal
basis of $\lql$, see [11, par 4].}\qed\medskip\noindent
{\bf 1.7.} Recall, [10, 6.4.27] that for all $\lambda$ in $P^+$,
the Kashiwara crystal $B(\lambda)$ is isomorphic to the Littelmann path
crystal
$C(\lambda)$. Moreover, in this isomorphism, the tensor product can be
defined in terms of
concatenation of paths [15].
Through this isomorphism, the following propositions are direct consequences
of [16, Theorem 10.1].
\medskip\noindent
{\bf Proposition A.} {\it Let $\lambda_i$ in $P^+$, $y_i$ in $W$,
 $b_i$ in $B(\lambda_i)$, with weight $y_i\lambda_i$, $1\leq i\leq r$.
 Suppose $y_1\leq y_2\leq\ldots\leq y_r$ for the Bruhat ordering. Then,
 $b_r\otimes b_{r-1}\otimes\ldots\otimes b_1$ belongs to the component of
type
 $B(\lambda_r+\lambda_{r-1}+\ldots+\lambda_1)$ in the crystal
 $B(\lambda_r)\otimes B(\lambda_{r-1})\otimes\ldots\otimes
B(\lambda_1)$.}\qed
 \medskip\noindent
 {\bf Proposition B.} {\it Let $\lambda$, $\mu$ in $P^+$,
$b$ in $B(\lambda)$, $c$ in $B(\mu)$ with weight $w_0\mu$.  Then,
 $c\otimes b$ belongs to the component of type
 $B(\mu+\lambda)$ in the crystal
 $B(\mu)\otimes B(\lambda)$.}\qed
 \medskip\noindent
For all $\lambda$ in $P^+$ and $w$ in $W$, let $B_w(\lambda)$ the subset of
$B(\lambda)$ corresponding to the Demazure module associated to $w$, see
[11, 12.4].\medskip\noindent
{\bf Proposition C.} {\it Let $\lambda$, $\mu$ in $P^+$, $w$ in $W$,
 $b$ in $B_w(\lambda)$, $c$ in $B_w(\mu)$ with weight $w\mu$.
  Then,
 $c\otimes b$ belongs to the component of type
 $B(\mu+\lambda)$ in the crystal
 $B(\mu)\otimes B(\lambda)$.}\qed
 \medskip\noindent
{\bf 1.8.} In this section, we define the string parametrization of the
canonical basis. For $b$ in $\B$ and $i$, $1\leq i\leq n$, set
$\tilde e_i^{max}(b)=\tilde e_i^{n_i(b)}(b)$, where
$n_i(b):=\hbox{max}\{k,\,\tilde e_i^{k}(b)\not=0\}$. Let $\tilde
w_0=s_{i_1}\ldots s_{i_N}$ be a reduced decomposition of $w_0$. To an
element $b$ in $\B$, we associate the element $A_{\tilde w_0}(b)$
in $\N^N$ such that the $k$-th component of $A_{\tilde w_0}(b)$ is
$n_{i_k}(\tilde e_{i_{k-1}}^{max}(\ldots(\tilde e_{i_1}^{max}(b)))$.
It is known, see [1], that $A_{\tilde w_0}$ is injective and that
$A_{\tilde w_0}(\B)$ is the set of integral points of a cone in $\R^N$.
Hence, the map $A_{\tilde w_0}$ defines a parametrization of the canonical
basis. It is called the string parametrization of $\B$.  \medskip\noindent
{\bf Proposition.} {\it We have
\smallskip\noindent
\item{(i)} Fix a reduced decomposition $\tilde w_0$ of the
longest
element of the Weyl group. Let $b$, $b'$, $b''$ be in $\B$ with
$A_{\tilde w_0}(b)+A_{\tilde w_0}(b')=A_{\tilde w_0}(b'')$.
Then, the $b''^*$-component of $b^*b'^*$ in the dual canonical
basis $\B^*$ is a power of $q$. \par\noindent
\item{(ii)} In the framework of
Propositions 1.7, there exists a reduced decomposition of $w_0$
for which the string parametrization of the tensor product is the sum of the
string parametrizations
of the factors.\smallskip\noindent
Proof.} (i) is a part of [8, Theorem 2.3].
Let's prove (ii). It is enough to prove it in the framework of Proposition
1.7. C. Indeed, the case A is obtained by induction from the case C and the
case B is obtained from the case C by setting $w=w_0$.
First suppose that $c$ and $b$ are respectively in $B(\mu)$ and
$B(\lambda)$ and that for some $i$, $\tilde f_i(c)=0$. Then , by [11, 4.3],
$$\varepsilon_i(c\otimes b)=\varepsilon_i(c)+\varepsilon_i(b)\leqno (*)$$
and $$\tilde e_i^{max}(c\otimes b)=\tilde e_i^{max}(c)\otimes\tilde
e_i^{max}(b).\leqno (**)$$ Let
$\tilde w=s_{i_1}\ldots s_{i_k}$ be a reduced decomposition of $w$ and let
$\tilde w_0=s_{i_1}\ldots s_{i_k}s_{i_{k+1}}\ldots s_{i_N}$ be the reduced
decomposition of $w_0$ which completes $\tilde w$ on the right hand.
For all $w'$ in $W$, let $c_{w'}$ be the extremal element of $B(\mu)$
corresponding to $w'$. By [11, 12.4], we know that $\tilde f_{i_1}(c_w)=0$ and
$\tilde e_{i_1}^{max}(c_w)=c_{s_{i_2}\ldots s_{i_k}}$.
Using (**) by induction, we obtain that
$\tilde e_{i_{k}}^{max}(\ldots(\tilde e_{i_1}^{max}(c_w\otimes b)))
=\tilde e_{i_{k}}^{max}(\ldots(\tilde e_{i_1}^{max}(c_w)))
\otimes\tilde e_{i_{k}}^{max}(\ldots(\tilde e_{i_1}^{max}(b)))$.
By [11, Theorem 12.4], the hypothesis $b\in B_w(\lambda)$ implies that
$\tilde e_{i_{k}}^{max}(\ldots(\tilde e_{i_1}^{max}(b)))$ is the highest
weight vector
$b_\lambda$. In particular, $\tilde e_{i_{k}}^{max}(\ldots(\tilde
e_{i_1}^{max}(c_w)))$ is the highest weight vector $b_\mu$.
Using (*) by induction, we obtain the proposition.
  \qed\medskip\noindent
  {\bf 1.9.} The following proposition is a given in [13, Proposition
33]. In [{\it loc. cit.}], the authors study the A$_n$ case, but the generalization
of their proof for all $\g$ is straightforward.\medskip\noindent
{\bf Proposition.} {\it Fix $r$ in $\N$ and $\lambda_i$,  $1\leq i\leq r$,
in $P^+$.
Let $b_i$, $1\leq i\leq r$, and $b$ respectively in $\B(\lambda_i)$,
$\B(\sum_i\lambda_i)$, such that \par
\item{(i)} $[b_1]\otimes[b_2]\ldots\otimes[b_r]$ is in the connected
component
of the crystal graph of $B(\lambda_1)\otimes
B(\lambda_2)\otimes\ldots\otimes
B(\lambda_r)$
of type  $B(\sum_i\lambda_i)$,
\par\noindent
\item{(ii)} $[b]$ identifies with
$[b_1]\otimes[b_2]\otimes\ldots\otimes[b_r]$ in
$B(\sum_i\lambda_i)$,\par\noindent
then there exists an integer $m$ such that the multiplication rule holds
in $\uqn$
$$q^mb_1^*b_2^*\ldots b_r^*=b^* \;\;\hbox{mod}\;\;q\L^*.$$}
\medskip\noindent
{\bf Corollary.} {\it Let $b_i$, $1\leq i\leq r$, $b$, $m$ as in the
previous
proposition. Moreover, suppose that \smallskip\noindent
\item{(iii)} the $b_i$ pairwise $q$-commute\par\noindent
\item{(iv)} for some $\tilde w_0$, the string parametrization of $[b]$ is
the sum of the string parametrizations
of the $[b_i]$. \smallskip\noindent
Then,
$q^mb_1^*b_2^*\ldots b_r^*=b^*$.\smallskip\noindent
Proof.} From the previous proposition, we obtain
 $$q^mb_1^*b_2^*\ldots b_r^*=b^*+q\sum_{b'^*\in\B^*}c_{b'^*}(q)b'^*,
 \leqno (1.9.1)$$
 where $c_{b'^*}(q)\in\Z[q]$. Clearly, the elements of the dual canonical
basis
 of the right hand term of (1.9.1) have the same weight equal to
$\sum_i\wt(b_i^*)$.
 Applying the automorphism $\eta$ and the antiautomorphism $\sigma$
on
 (1.9.1), Proposition 1.5 and (iii) give :
 $$q^mb_1^*\ldots b_{r-1}^*
b_r^*=q^kb^*+q^{k-1}\sum_{b'^*\in\B^*}c_{b'^*}(q^{-1})b'^*,$$
 for some integer $k$. Now, using (iv) and Proposition 1.8 (i),
we obtain that the sum in the right hand side of (1.9.1) does not contain
any term in $b^*$. Hence, by comparison, $k=0$ and then we
 obtain $c_{b'^*}(q)=0$. This gives the corollary.\qed\medskip\noindent
 Remark that, by a weight argument, the condition (iv) of the corollary
 implies (but is not equivalent to) the condition (i) of the proposition.\bigskip\noindent
{\bf 2. Adapted algebras} \bigskip\noindent
{\bf 2.1.} For each reduced decomposition of $w_0$, we define a set
of pairwise $q$-commuting elements of the dual canonical
basis.\medskip\noindent
{\bf Lemma.} {\it Fix a reduced decomposition $\tilde w_0
=s_{i_1}\ldots s_{i_N}$ of $w_0$. For all $r$, $1\leq r\leq N$, set
$y_r= s_{i_1}\ldots s_{i_r}$. Then,\par\noindent
\item{(i)} there exists a unique $X_{\tilde w_0}^r$ in $\uqn$ such
that    $v_{y_r\varpi_{i_r}}^*(uv_{\varpi_{i_r}})=
(K_{-\varpi_{i_r}}X_{\tilde w_0}^r,u)$,
for all $u$ in $\uqbm$,\par\noindent
\item{(ii)} up to a multiplicative scalar
$X_{\tilde w_0}^r$ is in $\B^*$,\par\noindent
\item{(iii)} the $X_{\tilde w_0}^r$, $1\leq r\leq N$, pairwise
$q$-commute.
\smallskip\noindent
Proof.} the existence of the q-commuting family of elements $\{X_{\tilde
w_0}^r
, \; 1\leq r\leq N\}$
is proved in [6, Proposition 3.2]. By [{\it loc. cit.}, Theorem 1.6],
$K_{-2\varpi_{i_r}}X_{\tilde w_0}^r$ is in ad$\uqn.K_{-2\varpi_{i_r}}$.
Hence, by [7, Theorem 1.6], $X_{\tilde w_0}^r$ is in the space generated
by $\B(\varpi_{i_r})^*$. By the Weyl character formula, see 1.3,
$X_{\tilde w_0}^r$ is equal to an element of  $\B(\varpi_{i_r})^*$ up
to a multiplicative scalar. This gives the lemma.\qed
\medskip\noindent
{\bf Remark and Definition.} We can choose the extremal vectors
$v_{y_r\varpi_{i_r}}$ such that $X_{\tilde w_0}^r\in \B^*$.
Let $\m^r$, such that $B(\m^r)^*=X_{\tilde w_0}^r$, $1\leq r\leq N$.
The elements $\m^r$ can be easily calculated at least
when the reduced decomposition $\tilde w_0$ corresponds to some quiver
orientation of
the Dynkin diagram of $\g$. Let $\m_k^r$ be the $k$-th entry of $\m^r$. By
[6, (3.2.2)] and (1.5.1)
:
$$ \m_k^r=\cases{1 & if $i_k=i_r$ and $k\leq r$\cr
0 & if not\cr}$$
It is likely that this formula generalizes for all decomposition
$\tilde w_0$.
\medskip\noindent
Let $A_{\tilde w_0}$ be the algebra generated by the
$X_{\tilde w_0}^r$, $1\leq r\leq N$.
Then, by [6, Corollary 3.2]\medskip\noindent
{\bf Proposition.} {\it Fix a reduced decomposition ${\tilde w_0}$ of
$w_0$. The algebra $A_{\tilde w_0}$ is an algebra of regular functions
on a quantum space. Moreover, let $S_{\tilde w_0}$ be the multiplicative
part
generated by the $X_{\tilde w_0}^r$, $1\leq r\leq N$.
Then $S_{\tilde w_0}$ is an Ore
set in
$A_{\tilde w_0}$ and $S_{\tilde w_0}^{-1}A_{\tilde w_0}=S_{\tilde
w_0}^{-1}\uqn$.}
\qed\medskip\noindent
{\bf 2.2.} Let's start with a definition.\medskip\noindent
{\bf Definition.} Fix a reduced decomposition $\tilde w_0$ of $w_0$. A
subalgebra
$A$ of $\C_q[U^-]$ will be called adapted if it is spanned by a subset
$P^*$ of $\B^*$ which verifies the
following properties :
\smallskip\noindent
\item{(i)} (multiplicativity) For all $p_i^*$ in $P^*$, $1\leq i\leq k$,
and all $(n_i)\in\N^k$, we have $\prod_i (p_i^*)^{n_i}\in\B^*$
up to a power of $q$,\par\noindent
\item{(ii)} (maximality) for all $q^*$ in  $\B^*\backslash P^*$,
there exists $p^*$ in $P^*$ such that
$p^*$ and $q^*$ are not multiplicative.\smallskip\noindent
We now prove the main theorem of this article. It is a generalization of
[21, Theorem 6.1] for all semi-simple Lie algebra $\g$ and for all
reduced decomposition of $w_0$.\medskip\noindent
{\bf Theorem.} {\it Fix a reduced decompositions $\tilde w_0$ of $w_0$.
The algebra $A_{\tilde w_0}$ is adapted. Moreover, the Ore set
$S_{\tilde w_0}$, see Proposition 2.1, is a subset of $\B^*$
which spans $A_{\tilde w_0}$.
\smallskip\noindent
Proof.} Let's prove that the algebra $A_{\tilde w_0}$ is adapted.
We know from Lemma 2.1 (iii) that the $X_{\tilde w_0}^r$ $q$-commute.
Hence, the
multiplicativity part is a consequence of Corollary 1.9. Indeed, the conditions of
the corollary are satisfied by Proposition A in
1.7 and Proposition 1.8 (ii). The last assertion of the theorem holds.\par
We now prove the maximality property.
Fix a reduced decomposition $\tilde w'_0$ which is compatible with a
quiver orientation. Let ${\cal C}_{\tilde w_0}^{\tilde w'_0}$ be the
set of parametrizations of elements of $S_{\tilde w_0}$ for the
Lusztig's parametrization associated to ${\tilde w'_0}$. By
Theorem 1.4, (1.5.1) and Proposition 2.1. ${\cal C}_{\tilde w_0}^{\tilde
w'_0}$
is the set of integral points of a simplicial cone in $\R^N$.
Suppose $q^*$ in $\B^*\backslash S_{\tilde w_0}$.
Then, by Proposition 2.1, there exists an
element $p^*$ in $S_{\tilde w_0}$ such that
$u:=p^*q^*\in A_{\tilde w_0}$. Suppose that $u$ is a monomial,
i.e. belongs to $S_{\tilde w_0}$ up to a multiplicative
scalar. Then, Theorem 1.4 and (1.5.1) imply
$\m+{\bf n}\in{\cal C}_{\tilde w_0}^{\tilde w'_0}$, where
$\m$ and ${\bf n}$ are respectively the ${\tilde w'_0}$-parametrizations
of $p^*$ and $q^*$. As $\m$ is in the simplicial cone
${\cal C}_{\tilde w_0}^{\tilde w'_0}$, we have ${\bf n}\in
{\cal C}_{\tilde w_0}^{\tilde w'_0}$, which contradicts the hypothesis
$q^*$ in $\B^*\backslash S_{\tilde w_0}$.
This ends the proof of the
Theorem.\qed
\medskip\noindent
In view of Proposition 2.1 and the previous theorem, we can see that the
Berenstein-Zelevinsky conjecture is true "up to localization".
\medskip\noindent
{\bf Remark.} By
[19, 14.2.5 and Lemma 1.2.8 (b)], the dual canonical basis is stable by
$\sigma$. Hence, for all adapted algebra $A$, the algebra $\sigma(A)$ is
adapted.
In the section 4.2, we give some examples for $\g$ of type A$_2$ and B$_2$.
\bigskip\noindent
{\bf 3. The multiplicativity property of the $q$-center.}\bigskip\noindent
{\bf 3.1.} Let $Z_q$ be the $q$-center of $\uqn$, i.e. the subspace of
$\uqn$
  generated
by elements which $q$-commute with all homogeneous elements of $\uqn$.
Let's describe $Z_q$ as a space and as an algebra. \medskip\noindent
{\bf Lemma.} {\it We have : \par\noindent
\item{(i)} For all $\lambda$ in $P^+$, there exists
 an element $z_\lambda$ in $\uqn$ such that
 $v_{w_0\lambda}^*(uv_{\lambda})=
(K_{-\lambda}z_\lambda,u)$,
for all $u$ in $\uqbm$,\par\noindent
\item{(ii)} $z_\lambda$ $q$-commutes with all homogeneous
elements of $\uqn$ and $Z_q$ is generated as a space by $z_\lambda$,
$\lambda\in P^+$,\par\noindent
\item{(iii)} $z_\lambda z_\mu=z_{\lambda+\mu}$ and so,
the algebra $Z_q$ is generated by the $z_{\varpi_k}$, $1\leq
k\leq n$,\par\noindent
\item{(iv)} up to a multiplicative scalar
$z_\lambda$ belongs to the dual canonical basis.\smallskip\noindent
Proof.} (i), (ii), and (iii) are proved in [4] and the proof of
(iv) is similar to the proof of Lemma 2.1 by remarking that $z_\lambda$
corresponds to the extremal vector $v_{w_0\lambda}^*$.\qed
\medskip\noindent
In the sequel, we choose $v_{w_0\lambda}$ such that $z_\lambda$
is in
$\B^*$. Fix a reduced decomposition $\tilde w_0=s_{i_1}\ldots s_{i_N}$
of $w_0$ and for all $k$, $1\leq k\leq n$, let $r(k)$ be the greatest
integer
such that $i_{r(k)}=k$. Then, it is easily seen that $z_{\varpi_k}=
X_{\tilde w_0}^{r(k)}$. Set ${\bf n}_k=\m_{r(k)}$. Then, when $\tilde w_0$
corresponds to a
quiver orientation, Remark 2.1 gives :
$$z_\lambda=B_{\tilde w_0}^*(\sum_k\lambda_k {\bf n}_k),\;\hbox{ where }\;
\lambda=\sum\lambda_k\varpi_k\in P^+.\leqno (3.1.1)$$
\medskip\noindent
{\bf 3.2.} As in the proof of Theorem 2.2, the following Proposition is a
consequence of Corollary 1.9, Proposition B in 1.7 and Proposition 1.8 (ii).\medskip\noindent
{\bf Proposition.} {\it Let $b^*$ be in $\B^*$ and let $c^*$ be in
$Z_q\cap\B^*$, then $b^*$ and $c^*$ are
multiplicative.}\qed\medskip\noindent
Let $J_q$ be the ideal of $\uqn$ generated by the $z_{\varpi_k}$, $1\leq
k\leq n$. The previous proposition implies that $J_q$ is generated as a
space by a part $\B^*(J_q)$ of $\B^*$. Let $H_q$ be the subspace generated
by the complementary set $\B^*(H_q)=\B^*\backslash
\B^*(J_q)$.\medskip\noindent
{\bf Corollary} {\it The $Z_q$-module $\uqn$ is  free with (canonical) basis
$\B^*(H_q)$, i.e. $\uqn=Z_q\otimes H_q$.
 The Berenstein-Zelevinsky conjecture is true on $\B^*$ if and only if it
is true on $\B^*(H_q)$.\smallskip\noindent
Proof.} Fix a reduced decomposition of $\tilde w_0=s_{i_1}\ldots s_{i_N}$
which is compatible with a quiver orientation.
 For all $k$, $1\leq k\leq n$, set
$I_k:=\{r,\;i_r=k\}$. remark that $\{1,\ldots,N\}$
is the disjoint union of the $I_k$, $1\leq i\leq n$.
Let $P_{\tilde w_0}(Z_q)$, resp. $P_{\tilde
w_0}(H_q)$, be the set of
parametrizations of $Z_q\cap\B^*$, resp. $\B^*(H_q)$, for $\tilde w_0$.
By Proposition 3.2, (3.1.1) and Remark 1.5,
it is easy to describe $P_{\tilde w_0}(Z_q)$ and
$P_{\tilde w_0}(H_q)$
:\smallskip\noindent
The semigroup $P_{\tilde w_0}(Z_q)$ is generated by the ${\bf n}_k$, $1\leq
k\leq n$,
whose
$r$-th entry is one if $r\in I_k$ and 0 if not.\par\noindent
The set $P_{\tilde w_0}(H_q)$ is the set of elements $\m$ of $\N^N$ such
that for all
$k$, $1\leq k\leq n$, there exists $r$ in $I_k$ such that the $r$-th entry
of $\m$ is 0.\smallskip\noindent
Now, for all $\m$ in $\N^N$, and $k$, $1\leq k\leq n$, let
$\m(k)$ be the minimal $r$-th entry of $\m$ when $r$ runs over $I_k$.
It is clear that $\m=(\sum_k\m(k){\bf n}_k)+(\m-\sum_k\m(k){\bf n}_k)$ is
the unique decomposition of $\m$ into a sum of an element of $P_{\tilde
w_0}(Z_q)$ and an
element of $P_{\tilde w_0}(H_q)$. This fact and the proposition implies the
freeness
assertion.\par
Now, suppose that the Berenstein-Zelevinsky conjecture is true on
$\B^*(H_q)$. Suppose that two elements $b^*$ and $b'^*$ of $\B^*$
$q$-commute. Let's prove that they are multiplicative. By the previous
proposition and the previous assertion, we can decompose $b^*=b_z^*b_h^*$,
 $b'^*=b_z'^*b_h'^*$, up to a power of $q$, where $b_z^*,\;b_z'^*\in
Z_q\cap\B^*$, $b_h^*,\;b_h'^*\in\B^*(H_q)$. By Theorem 1.4, two elements
in $\uqn$ commute up to a scalar imply that they $q$-commute.
If we suppose that $b_h^*$ and $b_h'^*$ do not $q$-commute, then we obtain
that there exists $u\not\in\K b_h^*b_h'^*$  such that
$b_h^*b_h'^*=q^mb_h'^*b_h^*+u$. This implies that
$b^*b'^*=q^{m_1}b'^*b^*+q^{m_2}ub_z^*b_z'^*$, and $ub_z^*b_z'^*\not\in\K
b^*b'^*$. This is in contradiction with the fact that $b^*$ and $b'^*$
$q$-commute. Hence, $b_h^*$ and $b_h'^*$  $q$-commute. By the hypothesis,
this implies that $b_h^*b_h'^*$ belongs to $\B^*$ up to a power of $q$.
Hence,
this also holds for $b_h^*b_h'^*b_z^*b_z'^*$ by the previous proposition.
This implies that $b^*$ and $b'^*$ are multiplicative.\qed\medskip\noindent
{\bf 3.3.} In this section, we prove that the intersection of all adapted
algebras
is  $Z_q$. To be more precise :\medskip\noindent
{\bf Proposition.} {\it Each adapted algebra contains $Z_q$. Moreover,
$\bigcap_{\tilde w_0}A_{\tilde w_0}=Z_q$,
where $\tilde w_0$ runs over the set of reduced decomposition of
$w_0$.\smallskip\noindent
Proof.} By Proposition 3.2 and by the maximality condition, each adapted
algebra
contains $Z_q$.
\par Remark that $\bigcap_{\tilde w_0}A_{\tilde w_0}$ is spanned by its
intersection with $\B^*$. Now, let $b^*\in
\bigcap_{\tilde w_0}A_{\tilde w_0}\cap\B^*$.
Then, by construction, $b^*$ $q$-commutes with
the
$X_{\tilde w_0}^r$. In particular, $b^*$ $q$-commutes with
$X_{\tilde w_0}^1=E_{i_1}$, where
$\tilde w_0=s_{i_1}s_{i_2}\ldots s_{i_N}$. Since each $s_k$, $1\leq
k\leq  n$
can be the first factor of a reduced decomposition of $w_0$, it implies that
$b^*$ is in $Z_q$. So, the proposition holds.\qed\medskip\noindent
Remark that in the intersection formula, we can choose a  subset
of cardinal n of the set $\tilde w_0$ of reduced decomposition of $w_0$.
\bigskip\noindent
{\bf 4. Adapted semigroups and examples.}\bigskip\noindent
{\bf 4.1.} Fix a reduced decomposition $\tilde w_0$ of $w_0$
which corresponds to a quiver orientation.
In this section, we study the set of parametrizations of
$A\cap\B^*$, where $A$ is an adapted algebra, for $\tilde
w_0$.\medskip\noindent
{\bf Proposition.} {\it Let $A$ be a adapted algebra, then the set
${\cal C}_{\tilde w_0}^A$ of parametrizations of $A\cap\B^*$ for
$\tilde w_0$ is a subsemigroup of $\N^N$.\smallskip\noindent
Proof.} This is a direct consequence of Theorem 1.4 and
(1.5.1).\qed\medskip\noindent
In the sequel, ${\cal C}_{\tilde w_0}^A$ will be called
adapted semigroup for ${\tilde w_0}$.
For any reduced decomposition $\tilde w'_0$ of $w_0$, set
${\cal C}_{\tilde w_0}^{\tilde w'_0}:={\cal C}_{\tilde w_0}^{A_{\tilde
w'_0}}$.
Proposition 2.1 and Theorem 2.2 give :\medskip\noindent
{\bf Proposition.} {\it The subset ${\cal C}_{\tilde w_0}^{\tilde w'_0}$
of $\N^N$
is the set of integral points of a simplicial cone of $\R^N$.
Moreover, ${\cal C}_{\tilde w_0}^{\tilde w'_0}$ generates $\Z^N$ as a
group.}\qed\medskip\noindent
We do not know if an adapted semigroup is in general
the set of integral points of a real cone. If this case occurs,
 it will be called adapted cone.\medskip\noindent
{\bf 4.2.} In this section, we study the cases where $\g$ is of type A$_2$
and B$_2$. For each case, we fix a reduced decomposition ${\tilde w_0}$ for
the parametrization of the dual canonical basis and then, we give explicitly
the generators of the cones $\C_{\tilde w_0}^{\tilde w'_0}$ for
all ${\tilde w'_0}$. The generators of $\C_{\tilde w'_0}^{\tilde w'_0}$ are
calculated from Remark 2.1 and the "reparametrization" $\C_{\tilde
w_0}^{\tilde w'_0}$ is calculated with the help of [1, Proposition 7.1 (3)].
The images
$\sigma(\C_{\tilde w_0}^{\tilde w'_0})$ can be easily guessed by using
[11, Proposition 8.2].\medskip\noindent
{\bf Example 1.} Let $\g$ of type A$_2$. Fix ${\tilde w_0}=s_1s_2s_1$,
 ${\tilde w'_0}=s_2s_1s_2$.
Then, $\C_{\tilde w_0}^{\tilde w_0}$ is generated by $\{(1,0,0),\,
(0,1,0),\, (1,0,1)\}$. $\C_{\tilde w'_0}^{\tilde w_0}$ is generated by
$\{(0,0,1),\,
(0,1,0),\, (1,0,1)\}$. Both cones are stable by $\sigma$ and the union
$\C_{\tilde w_0}^{\tilde w_0}\cup\C_{\tilde w'_0}^{\tilde w_0}$ is the whole
semigroup $\N^N$.\medskip\noindent
{\bf Example 2.} Let $\g$ of type B$_2$ with minuscule weight $\varpi_2$.
Fix ${\tilde w_0}=s_2s_1s_2s_1$,
 ${\tilde w'_0}=s_1s_2s_1s_2$.
 \par\noindent
 The previous construction gives four cones, with intersection
 $P_{\tilde w_0}(Z_q)$ generated by $(1,0,1,0)$ and $(0,1,0,1)$.
 In the following, we present the four cones  $\C_{\tilde w_0}^{\tilde
w_0}$,
 $\C_{\tilde w'_0}^{\tilde w_0}$,
 $\sigma(\C_{\tilde w_0}^{\tilde w_0})$,
 $\sigma(\C_{\tilde w'_0}^{\tilde w_0})$ and their canonical set of
generators :
 $$\C_{\tilde w_0}^{\tilde w_0}\hskip 1cm\{(1,0,0,0),\,
(0,1,0,0),\, (1,0,1,0),\,(0,1,0,1)\},$$
$$\C_{\tilde w_0}^{\tilde w'_0}\hskip 1cm\{(0,0,0,1),\,
(1,0,0,1),\, (1,0,1,0),\,(0,1,0,1)\},$$ $$\sigma(\C_{\tilde w_0}^{\tilde
w_0})
\hskip 1cm\{(1,0,0,0),\,
(2,0,0,1),\, (1,0,1,0),\,(0,1,0,1)\},$$ $$\sigma(\C_{\tilde w_0}^{\tilde
w'_0})
\hskip 1cm\{(0,0,0,1),\,
(0,0,1,0),\, (1,0,1,0),\,(0,1,0,1)\},$$
where $\sigma$ is the map which commutes with the parametrization
corresponding
with $\tilde w_0$.
With the same methods used in 2.1 and 2.2, we obtain
that the semigroup ${\cal D}$ generated by
$${\cal D}\hskip 1cm\{(0,1,0,0),\,
(0,0,1,0),\, (1,0,1,0),\,(0,1,0,1)\}$$ is an adapted cone. Hence,
$\sigma({\cal D})$, which is generated by
$$\sigma({\cal D})\hskip 1cm\{(2,0,0,1),\,
(1,0,0,1),\, (1,0,1,0),\,(0,1,0,1)\},$$ is also
an adapted cone. Let's give a few details.
\medskip\noindent
{\bf Lemma.} {\it The cones ${\cal D}$ and $\sigma({\cal D})$
are adapted.\smallskip\noindent
sketch of the proof.} It is enough to prove that
${\cal D}$ is an adapted cone, see Remark 2.2. Set
${\bf a}=(0,0,1,0)$, ${\bf b}=(0,1,0,0)$,
${\bf c}=(1,0,1,0)$, ${\bf d}=(0,1,0,1)$. By (3.1.1),
$B_{\tilde w_0}^*({\bf c})$ $B_{\tilde w_0}^*({\bf d})$
generate the $q$-center $Z_q$. By 3.2, in order to prove the
multiplicativity
part, it is enough to
prove that, up to a power of $q$,
$B_{\tilde w_0}^*({\bf a})^aB_{\tilde w_0}^*({\bf b})^b$
is in $\B^*$ for all non negative integer $a$ and $b$. Through
Corollary 1.9 and Proposition 1.7. C, this is implied by the
following assertions :
\smallskip\noindent
1) $B_{\tilde w_0}^*({\bf a})$ and $B_{\tilde w_0}^*({\bf b})$ q-commute,
\par\noindent
2) $B_{\tilde w_0}^*({\bf a})^a$ is in $\B^*$ up to a power of
$q$,\par\noindent
3) In the $W$-stratification of $\B^*$, see [11, 12.4], $B_{\tilde w_0}^*({\bf
a})^a$
belongs to the $s_2s_1$ component of $\B^*$ and
$B_{\tilde w_0}^*({\bf b})$ corresponds to the extremal vector
 $v_{s_2s_1\varpi_1}^*$ of $B(\varpi_1)^*$.\smallskip\noindent
 Let sketch the proof for 1) 2) and 3).
 By Remark 1.5, we can describe the parametrizations
of $B(\varpi_i)^*$ for $i=1,2$. We then obtain that
$B_{\tilde w_0}^*({\bf a})$ and $B_{\tilde w_0}^*({\bf b})$
are in $B(\varpi_1)^*$.
 Let $v^*$ and $w^*$ be weight
 vectors in $V_q(\varpi_1)^*$ corresponding respectively to
 $B_{\tilde w_0}^*({\bf a})$ and $B_{\tilde w_0}^*({\bf b})$.
 It is easy to prove that for all positive root $\alpha$, then
 $E_\alpha.v^*$ is non zero implies that $F_\alpha.w^*$ is zero.
 By [5, 1.5], this implies 1). Now, by [11, Proposition 8.2], $\sigma(
B_{\tilde w_0}^*({\bf a}))$ is in $B(\varpi_2)^*$, so it corresponds to an
extremal
vector. By Proposition 1.7. A, Proposition 1.8 (ii)
 and Corollary 1.9, it gives that
$\sigma(
B_{\tilde w_0}^*({\bf a}))^a$ is in $\B^*$. This implies 2).
Moreover, we can obtain by considering the crystal graph of
$V_q(\varpi_1)$ that $B_{\tilde w_0}^*({\bf a})$
belongs to the $s_2s_1$ component of $\B^*$. This implies 3).
\par\noindent
The maximality property is proved as in Theorem 2.2 by remarking that
${\cal D}$ generates $\Z^N$ as a group.
\qed\medskip\noindent
Let ${\cal G}^*$ be the finite subset of $\B^*$ parametrized
via ${\tilde w_0}$ by  the union of the canonical
generators of the cones 
$\C_{\tilde w_0}^{\tilde w_0}$,
$\C_{\tilde w_0}^{\tilde w'_0}$, $\sigma(\C_{\tilde w_0}^{\tilde w_0})$,
$\sigma(\C_{\tilde w_0}^{\tilde w'_0})$, ${\cal D}$, $\sigma({\cal
D})$, then  
 it is easily verified that \medskip\noindent
 {\bf Proposition.} {\it We have the decomposition
 $\N^N= \C_{\tilde w_0}^{\tilde w_0}\cup
\C_{\tilde w_0}^{\tilde w'_0}
 \cup\sigma(\C_{\tilde w_0}^{\tilde w_0})
 \cup\sigma(\C_{\tilde w_0}^{\tilde w'_0})\cup{\cal D}\cup \sigma({\cal
D})$. Hence, each element of $\B^*$ can be decomposed into a product of elements
of ${\cal G}^*$.}\qed
 \vskip 1 cm \noindent{\petittitre\centerline{ACKNOWLEDGMENTS}}\vskip 1cm
 We wish to thank Bernard Leclerc for useful conversations on the
conjecture.
 We are grateful to Olivier Mathieu for an interpretation of the condition (i) of
 Proposition 1.9 in terms
 of Littelmann's path model.
\vskip 1 cm \noindent{\petittitre\centerline{BIBLIOGRAPHY}}\vskip 1cm
\parindent=0cm
[1]. A. BERENSTEIN and A. ZELEVINSKY. {\it Tensor product
multiplicities, Canonical bases
and Totally positive varieties}, Invent. Math. 143 (2001), no. 1, 77--128.
\smallskip
[2]. A. BERENSTEIN and A. ZELEVINSKY. {\it String bases for quantum groups
of type $A\sb r$},
I. M. Gel'fand Seminar, 51--89,
Adv. Soviet Math., 16, Part 1,
Amer. Math. Soc., Providence, RI, (1993).\smallskip
[3]. P. CALDERO. {\it Elements ad-finis de certains groupes
quantiques},
 C.R.Acad.
Sci. Paris, t. 316, Serie I, (1993), 327-329.\smallskip
[4]. P. CALDERO. {\it Etude des $q$-commutations dans l'alg\`ebre
$U_q(\n)$,}
J. Algebra, 178, (1995), 444-457.\smallskip
[5]. P. CALDERO. {\it On the Gelfand-Kirillov conjecture
for quantum algebras}, Proc. Amer. Math. Soc., 128, (2000),
943--951.\smallskip
[6]. P. CALDERO. {\it On the $q$-commutations in $U_q(\n)$ at roots of
one},
J. Algebra,
210, (1998), 557-576.\smallskip
[7]. P. CALDERO. {\it On harmonic elements for semi-simple Lie algebras},
accepted in Adv. in Maths, math.RT/0012092.\smallskip
[8]. P. CALDERO. {\it Toric degenerations of Schubert varieties}, 
math.RT/0012165 .\smallskip
[9]. C. DE CONCINI and V. G. KAC. {\it Representations of quantum groups at
roots of 1}, Colloque Dixmier, Progress in Math., 92, (1990), 471-506.
\smallskip
[10]. A. JOSEPH. {\it Quantum groups and their primitive ideals,}
 Springer-Verlag, 29, Ergebnisse der Mathematik und ihrer Grenzgebiete,
(1995).
 \smallskip
[11]. M. KASHIWARA. {\it On Crystal Bases}, Canad. Math. Soc.,
Conference Proceed., 16, (1995), 155-195.\smallskip
[12] B. LECLERC and A. ZELEVINSKY. {\it Quasicommuting families of
quantum Plücker coordinates}, Kirillov's seminar on representation theory,
85--108,
Amer. Math. Soc. Transl. Ser. 2, 181, (1998).\smallskip
[13] B. LECLERC, M. NAZAROV and J-Y THIBON. {\it Induced
representations of affine Hecke algebras and canonical bases of quantum
groups},
ArXiv:Math.QA/0011074.\smallskip
[14]. S.Z. LEVENDORSKII, Y.S. SOIBELMAN. {\it Some applications of quantum
Weyl
 group,} J. Geom. Phys., 7, (1990), 241-254.\smallskip
[15]. P. LITTELMANN,
{\it A Littlewood-Richardson rule for
symmetrizable Kac-Moody algebras}, Invent. Math., 116, (1994),
329--346.\smallskip
[16]. P. LITTELMANN,
{\it A plactic algebra for
semisimple Lie algebras}, Adv. Math. 124 (1996), no. 2, 312--331.\smallskip
[17]. G. LUSZTIG, {\it Canonical bases arising from quantized enveloping
algebras},
 J. Amer. Math. Soc. 3 (1990), no. 2, 447--498.\smallskip
[18]. G. LUSZTIG. {\it Canonical bases arising from quantized enveloping
algebras.
 II}, Prog. Theor. Phys. Suppl., 102, (1990), 175-201.
\smallskip
[19].  G. LUSZTIG. {\it Introduction to quantum groups}, Progress in
Mathematics,
 110, Birkauser, (1993).
\smallskip
[20]. M. REINEKE. {\it On the coloured graph
structure of Lusztig's canonical basis}, Math. Ann., 307, (1997),
 705--723.\smallskip
[21]. M. REINEKE. {\it Multiplicative Properties of Dual Canonical Bases
of
Quantum Groups}, J. Alg., 211, (1999), 134-149.\smallskip
[22]. M. ROSSO. {\it Analogues de la forme de Killing et du th\'eor\`eme
 de Harish-Chandra
pour les groupes quantiques}, Ann. Sci. Ec. Norm. Sup., 23, (1990),
445-467.\smallskip
[23]. Y. SAITO. {\it PBW basis of quantized universal 
enveloping algebras}, Publ. Res. Inst. Math. Sci., 30, (1994), 
209--232.\smallskip 
[24]. T. TANISAKI. {\it Killing forms, Harish-Chandra isomorphisms, and
universal
{\cal R}-matrices for quantum algebras}, Int. J. Mod. Phys. A, Vol. 7,
 Suppl. 1B,
 (1992),  941-961.\smallskip
{\it Institut Girard Desargues, UPRES-A-5028  \par
Universit\'e Claude Bernard Lyon I, Bat 101\par
69622 Villeurbanne Cedex, France}\smallskip\noindent
e-mail : caldero@desargues.univ-lyon1.fr
\medskip\noindent
Keywords : Quantum groups, Canonical Basis.
\end